\documentclass[11pt,reqno]{amsart}
\usepackage{amsmath,amsfonts,amssymb,amsthm,stmaryrd,enumerate,url}
\usepackage{microtype}
\usepackage{bbm,mathrsfs}
\usepackage[utf8]{inputenc}  %% pour les accents
\usepackage[english]{babel}
\usepackage[numbers,sort]{natbib}
\usepackage[usenames,dvipsnames,svgnames,table]{xcolor}
\usepackage[colorlinks=true, pdfstartview=FitV, pagebackref=true]{hyperref}
\usepackage{enumitem}

\usepackage[a4paper,vmargin=4cm,inner=3cm,outer=3cm]{geometry}

\newtheorem{theorem}{Theorem}

\newtheorem{proposition}[theorem]{Proposition}

\newcommand{\R}{\mathbb R}

\newcommand{\Z}{\mathbb Z}

\newcommand{\E}{\mathbb E}

\newcommand{\cE}{{\mathcal E}}

\DeclareMathOperator{\var}{Var}

\DeclareMathOperator{\ent}{Ent}

\newcommand{\ind}{{\mathbf{1}}}
\newcommand{\tmix}{t_\mathsf{mix}}

\newcommand{\cls}{\mathrm{C}_\mathsf{LS}}
\newcommand{\tv}{\mathsf{tv}}

\begin{document}

\title[Mixing time of a matrix random walk]{Mixing time of a matrix random walk generated by elementary transvections}

\author[]{Anna Ben-Hamou}
\address[]{Sorbonne University, LPSM, \\
4 place Jussieu, 75005 Paris, France}
\email[]{anna.ben\_hamou@sorbonne-universite.fr}

\begin{abstract}
We consider a Markov chain on invertible~$n\times n$ matrices with entries in~$\Z_2$ which moves by picking an ordered pair of distinct rows and add the first one to the other, modulo~$2$. We establish a logarithmic Sobolev inequality with constant~$n^2$, which yields an upper bound of~$O(n^2\log n)$ on the mixing time.
\end{abstract}

\maketitle

\section{Setting and main result}

We consider a random walk on~$\Omega=\mathrm{SL}_n(\Z_2)$, the set of invertible binary matrices, generated by elementary transvections: at each time, the walk proceeds by picking two distinct rows at random, and adding one to the other, modulo~$2$. Letting~$E_{i,j}$ be the~$n\times n$ matrix with a~$1$ in position~$(i,j)$, and~$0$ elsewhere, the transition matrix~$P$ of the walk is defined by
\begin{displaymath}
P(x,y)\;=\;\begin{cases}
\frac{1}{n(n-1)} &\mbox{if $y=(I + E_{i,j})x$ for some $1\leq i\neq j\leq n$,}\\
0 &\mbox{otherwise.}
\end{cases}
\end{displaymath}

Since the elementary transvections~$I + E_{i,j}$ generate~$\Omega$, the chain is irreducible. By symmetry of the matrix~$P$, the uniform distribution~$\pi$ over~$\Omega$ is the unique stationary distribution. Let us point out here that~$\Omega$ is a significant subset of the set of all binary matrices. Indeed
\[
\frac{|\Omega|}{2^{n^2}}\;=\;\prod_{k=1}^n\left(1-2^{-k}\right) \;\approx\; 0.288788 \, .
\]
Since the chain is also aperiodic (see~\cite[Proposition 4.4]{diaconis1996walks}), the chain converges to~$\pi$ as time tends to infinity. To quantify the speed of this convergence, one usually considers the total variation mixing time defined, for all~$\varepsilon\in (0,1)$, by
\[
\tmix(\varepsilon)\;=\;\min\{t\geq 0\, ,\, \max_{x\in\Omega}\| P^t(x,\cdot)-\pi\|_\tv\leq \varepsilon\}\, ,
\]
where $\|\mu-\nu\|_\tv=\max_{A\subset \Omega}|\mu(A)-\nu(A)|$ is the total-variation distance. A stronger notion is the $\ell_2$-mixing time defined, for all~$\varepsilon\in (0,1)$, by
\[
\tmix^{(2)}(\varepsilon)\;=\;\min\{t\geq 0\, ,\, \max_{x\in\Omega}\left\| \frac{P^t(x,\cdot)}{\pi(\cdot)}-1\right\|_{2,\pi}\leq \varepsilon\}\, ,
\]
where $\|f\|_{2,\pi}=\left(\sum_{\in\Omega}\pi(x)f(x)^2\right)^{1/2}$ is the $\ell_2$-norm with respect to $\pi$. By Cauchy-Schwarz inequality, we have
\[
\tmix(\varepsilon)\;\leq\;  \tmix^{(2)}(2\varepsilon)\, .
\]

We obtain the following bound.

\begin{theorem}\label{thm:mixing}
There exists absolute constants $A>0$ and $B>0$ such that for all $\varepsilon>0$, the $\ell_2$-mixing time of the chain $P$ satisfies
\[
\tmix^{(2)}(\varepsilon)\;\leq \; A n^2\log n +B n \log\left(\frac{1}{\varepsilon}\right)\, .
\]
\end{theorem}

Theorem~\ref{thm:mixing} comes from 	an upper bound of $O(n^2)$ on the logarithmic Sobolev constant~$\cls$ of the chain, defined as the smallest constant~$c>0$ such that for all functions~$f\colon \Omega\to \R$, we have
\[
\ent_\pi(f^2)\;\leq\; c\,\cE_P(f, f)\, ,
\]
where~$\ent_\pi(f^2)=\E_\pi[f^2\log f^2] -\E_\pi[f^2]\log\E_\pi[f^2]$ is the entropy of~$f^2$, and where~$\cE_P(f,f)$ is the Dirichlet form defined by
\[
\cE_P(f,f) \;=\;\frac{1}{2}\sum_{x,y\in\Omega} \pi(x)P(x,y) \left(f(x)-f(y)\right)^2\, .
\]
By~\cite[Corollary 3.8]{diaconis1996logarithmic}, we have
\[
\tmix^{(2)}(\varepsilon)\;\leq\; \frac{\cls}{4}  \log\log\left(\frac{1}{\pi_\star}\right)+\lambda_*^{-1}\log\left(\frac{\sqrt{1+2e^2}}{\varepsilon}\right)+1\, ,
\]
where~$\pi_\star=\min_{x\in\Omega}\pi(x)$ and $\lambda_*$ is the absolute spectral gap of $P$. By~\cite[Proposition 4.4]{diaconis1996walks} and the powerful results of~\cite[Section 3.2]{kassabov2005kazhdan}, we have $\lambda_*^{-1}=O(n)$. Finally, since $\pi_\star\asymp \frac{1}{2^{n^2}}$, we see that Theorem~\ref{thm:mixing} will follow from the following proposition. 

\begin{proposition}\label{prop:log-sob}
The logarithmic Sobolev constant of $P$ satisfies
\[
\cls\; \lesssim\; n^2\, .
\]
\end{proposition}

In the above proposition and in the rest of the paper, the notation $u_n\lesssim v_n$ means that there exists an absolute constant $c>0$ such that for all $n\geq 1$, we have $u_n\leq cv_n$.

Interestingly, a crucial ingredient in the proof of Proposition~\ref{prop:log-sob} is the bound of $O(n)$ on the Poincar\'e constant (inverse of the spectral gap), obtained by~\cite{kassabov2005kazhdan}.

\section{Related work}

To our knowledge, the chain~$P$ was first studied by~\cite{diaconis1996walks}, who used a comparison with another chain on~$\Omega$ studied by~\cite{hildebrand1992generating} to get an upper bound of~$O(n^2)$ on the relaxation time (the inverse of the spectral gap), implying an upper bound of~$O(n^4)$ on the $\ell^2$-mixing time. Then the results of~\cite{kassabov2005kazhdan} on the Kazhdan constant yield an upper bound of~$O(n)$ on the relaxation time (inverse spectral gap), and thus of~$O(n^3)$ on the $\ell^2$-mixing time, which was the best known upper bound for total variation mixing as well. As for the lower bound, a simple counting argument shows that the total-variation mixing time can be lower bounded by~$\Omega\left(\tfrac{n^2}{\log n}\right)$ (which is actually an estimate of the diameter of the underlying graph, see \cite{andren2007complexity,christofides2014asymptotic}). It would be interesting to know whether a lower bound of~$\Omega(n^2\log n)$ holds, and more ambitiously, whether this chain exhibits the cutoff phenomenon.

This matrix random walk has also received interest from cryptographers, who have used it in authentication protocols. In cryptography, an authentication protocol involves two parties, a verifier and a prover. The verifier's goal is to confirm the prover's identity and to distinguish between an honest and a dishonest prover. A significant subset of authentication protocols, known as time-based authentication protocols, relies on the time it takes the prover to respond to a challenge. Authentication is successful only if the prover answers correctly and fast enough.

The protocol proposed by~\cite{sotiraki2016authentication} works as follows: starting with the identity matrix in~$\text{SL}_n(\Z_2)$, the prover runs the chain with kernel~$P$ for a specific duration~$t$. The resulting matrix~$A_t$ is made public (acting as the public key), but only the prover knows the chain's trajectory. Then, when he needs to authenticate, the prover requests a binary vector of length~$n$ from the verifier, and the challenge is to compute~$y=A_t x$ quickly. Since the prover is aware of the chain's trajectory, he can apply the same row operations to~$x$ that were used to generate~$A_t$, thus providing the correct answer in time~$t$. On the other hand, if~$t$ is sufficiently large, a dishonest party may struggle to distinguish~$A_t$ from a matrix chosen uniformly at random within polynomial time (indicating that the chain is \emph{computationally mixed}). In this case, the dishonest party's best option would be to perform regular matrix-vector multiplication, which typically requires~$n^2$ operations. Therefore, if the prover selects~$t$ as the computational mixing time, and if~$t$ is proven to be significantly smaller than~$n^2$, the verifier can successfully distinguish between honest and dishonest parties. Given the present result on the total variation mixing time $O(n^2 \log n)$, finding the computational mixing time (when only test functions that are computable in polynomial time are allowed) would be very interesting in this context.

One may notice that the restricted dynamics on the first~$k\leq n$ columns of the matrix is still a Markov chain. For~$k=1$, \cite{ben2018cutoff} established a cutoff at time~$\tfrac{3}{2}n\log n$ with a window of order~$n$. Later, the results of~\cite{peres2020cutoff} on finite groups allowed to generalize this cutoff for a fixed number~$k$ of columns. 

Another related chain is the random walk on upper triangular matrices with entries in~$\Z_q$ (where~$q\geq 2$ is an integer) and ones along the diagonal, where at each step a uniformly chosen row is added to the row above. For~$q=2$, \cite{peres2013mixing} showed that the mixing time is~$O(n^2)$, and \cite{nestoridi2023random} then obtained optimal bounds for the mixing time both in~$n$ and~$q$ (we refer the reader to this last paper for more reference on this process).

\section{Proof of Proposition~\ref{prop:log-sob}}

Our goal is to show that for all functions $f\colon \Omega\to \R$, we have
\[
\ent_\pi(f^2)\;\lesssim\; n^2\cE_P(f, f)\, .
\]
So let $f\colon\Omega\to \R$ and let $g$ be the extension of $f$ over the set $\Omega^*$ of all $n\times n$ matrices with entries in $\Z_2$ defined by
\[
\forall x\in \Omega^*\, ,\; g(x)\;=\;
\begin{cases}
f(x) &\text{ if $x\in\Omega$,}\\
\E_\pi[f] &\text{ if $x\not\in\Omega$.}
\end{cases}
\]
By the variational characterization of entropy, we have
\begin{eqnarray*}
\ent_\pi(f^2)&=& \inf_{u\in \R_+} \sum_{x\in \Omega}\frac{1}{|\Omega|}\left( f^2(x)\log\left(\frac{f^2(x)}{u}\right) -f^2(x) +u\right)\, .
\end{eqnarray*}
Using the fact that $t\log\left(\frac{t}{u}\right) -t+u\geq 0$ for all $t,u\in\R_+$, one may extend the sum to $\Omega^*$ and replace $f$ by $g$ (which coincide with $f$ on $\Omega$) without increasing the value. More precisely, we have
\begin{eqnarray*}
\ent_\pi(f^2) &\leq &  \inf_{u\in \R_+} \sum_{x\in \Omega^*}\frac{1}{|\Omega|} \left( g^2(x)\log\left(\frac{g^2(x)}{u}\right) -g^2(x) +u\right)\, .
\end{eqnarray*}
In particular, denoting by $\mu$ the uniform measure over $\Omega^*$ and taking $u=\E_\mu[g^2]$, we have
\begin{eqnarray*}
\ent_\pi(f^2)&\leq & \sum_{x\in \Omega^*}\frac{1}{|\Omega|}  \left( g^2(x)\log\left(\frac{g^2(x)}{\E_\mu[g^2]}\right) -g^2(x) +\E_\mu[g^2]\right)\\
&=  &\frac{2^{n^2}}{|\Omega|} \ent_\mu (g^2)\, .
\end{eqnarray*}
If $X$ is a random matrix with distribution~$\mu$, the rows $L_1,\dots,L_n$ of $X$ are independent, so by the sub-additivity of entropy for product measures, we have
\begin{eqnarray*}
\ent_\mu(g^2) &\leq & \sum_{i=1}^n \E_\mu\left[\ent_\mu\left( g^2\,|\, (L_k)_{k\neq i}\right)\right]
\end{eqnarray*}
For each $i\in [n]$, conditionally on $(L_k)_{k\neq i}=(\ell_k)_{k\neq i}$, observe that if the family $(\ell_k)_{k\neq i}$ does not span a subset of dimension $n-1$, then, whatever the value of $L_i$, the function $g$ is equal to $\E_\pi[f]$, so the conditional entropy is zero. Otherwise, if the family $(\ell_k)_{k\neq i}$ spans a subset of dimension $n-1$, then one may complete it with some vector $\ell_i$ to form a basis of $\{0,1\}^n$. The conditional law of $L_i$ is then that of $\sum_{k=1}^n a_k\ell_k$, where $(a_k)_{k=1}^n$ is uniformly distributed over $\{0,1\}^n$. Since the simple random walk on the hypercube has logarithmic Sobolev constant equal to~$n$ (see for instance~\cite[Example 3.2]{diaconis1996logarithmic}), we have, for each $i\in [n]$,
\[
\ent_\mu\left( g^2\,|\, (L_k)_{k\neq i}=(\ell_k)_{k\neq i}\right)\leq \frac{1}{2} \sum_{j=1}^n \E_\mu\left[\left(g(X)-g(X^{i\leftarrow j})\right)^2 \,\big|\, (L_k)_{k\neq i}=(\ell_k)_{k\neq i}\right]\, ,
\]
where $X^{i\leftarrow j}$ is the matrix obtained from $X$ by replacing $L_i$ by $L_i+\ell_j$. There are three different situations:
\begin{itemize}
\item either $X\in \Omega$ ($a_i=1$) and $j\neq i$: then $X^{i\leftarrow j}\in\Omega$ and we recover the Dirichlet form of $f$;
\item or $X\not\in \Omega$ ($a_i=0$) and $j\neq i$: then $X^{i\leftarrow j}\not\in\Omega$ and the difference is zero;
\item or $j=i$: then either $X\in\Omega$ and $X^{i\leftarrow j}\not\in\Omega$ or the other way round, and we recover the variance of $f$.
\end{itemize}
More precisely, we obtain
\[
\ent_\mu(g^2)\;\leq \; \frac{1}{2} \sum_{i=1}^n \sum_{j\neq i} \E_\mu\left[\ind_{X\in\Omega}\left(f(X)-f(X^{i\leftarrow j})\right)^2\right] +\sum_{i=1}^n \E_\mu\left[\ind_{X\in\Omega}\left(f(X)-\E_\pi[f]\right)^2\right]
\]
Coming back to the entropy of~$f^2$ under the measure $\pi$, this gives
\begin{eqnarray*}
\ent_\pi(f^2) &\leq & \frac{1}{2} \sum_{i=1}^n \sum_{j\neq i} \E_\pi\left[\left(f(X)-f(X^{i\leftarrow j})\right)^2\right] + n\E_\pi\left[\left(f(X)-\E_\pi[f]\right)^2\right]\\
& = & n(n-1)\cE_P(f,f)+n\var_\pi(f)\, .
\end{eqnarray*}
One may now conclude with the results of~\cite[Section 3.2]{kassabov2005kazhdan} on the spectral gap, which yield 
\[
\var_\pi(f)\;\leq\; 4(31\sqrt{n}+700)^2\cE_P(f,f)\, .
\]

\medskip

\noindent\textbf{Acknowledgments.} The author would like to thank Yuval Peres for suggesting this question back in 2015 during an internship under his supervision, and for numerous enlightening discussions.

\medskip

\noindent\textbf{Funding.} The author is supported by the French Agence Nationale de la Recherche, under the grants ANR-21-CE40-0019 (CORTIPOM project) and ANR-23-CE40-0003 (CONVIVIALITY project).

\bibliographystyle{abbrvnat}
\bibliography{biblio}

\end{document}